\documentclass{ceb-l}

\usepackage{graphicx}
\usepackage{subcaption} 
\usepackage{color, enumerate} 
\usepackage[multiple]{footmisc}
\newtheorem{theorem}{Theorem}[section]

\theoremstyle{definition}
\newtheorem{definition}[theorem]{Definition}

\newtheorem{question}[theorem]{Question}
\newtheorem{conjecture}[theorem]{Conjecture}

\theoremstyle{remark}

\numberwithin{equation}{section}

\newcommand{\C}{\mathbb{C}}

\newcommand{\N}{\mathbb{N}}

\newcommand{\R}{\mathbb{R}}

\newcommand{\T}{\mathbb{T}}

\newcommand{\cL}{\mathcal{L}}

\newcommand{\fkI}{\mathfrak{I}}

\newcommand{\ud}{\mathrm{d}}

\begin{document}

\title[The Kakeya Conjecture]{The Kakeya Conjecture: where does it come from and why is it important?}


%    author one information
\author{Jonathan Hickman}
\address{School of Mathematics and Maxwell Institute for Mathematical Sciences, James Clerk Maxwell Building, The King's Buildings, Peter Guthrie Tait Road, Edinburgh, EH9 3FD, UK.}
%\curraddr{}
\email{jonathan.hickman@ed.ac.uk}
\thanks{Supported by New Investigator Award UKRI097. This article expands upon the Oberwolfach snapshot \cite{Hickman}. The author thanks Tony Carbery, Tom Leinster, Betsy Stovall and Joshua Zahl for helpful comments and feedback on an earlier draft.}

%    author two information
%\author{}
%\address{}
%\curraddr{}
%\email{}
%\thanks{}

%    \subjclass is required.

\date{}

\dedicatory{}

%    Abstract is required.
\begin{abstract}
Roughly speaking, the Kakeya Conjecture asks to what extent lines which point in different directions can be packed together in a small space. In $\R^2$, the problem is relatively straightforward and was settled in the 1970s. In $\R^3$ it is much more difficult and was only recently resolved in a monumental and groundbreaking work of Hong Wang and Joshua Zahl. This note describes the origins of the Kakeya Conjecture, with a particular focus on its classical connections to Fourier analysis, and concludes with a discussion of elements of the Wang--Zahl proof. The goal is to give a sense of why the problem is considered so central to mathematical analysis, and thereby underscore the importance of the Wang--Zahl result. 
\end{abstract}

\maketitle


\bibliographystyle{amsplain}

\section{Kakeya sets}

\subsection{Kakeya's question} A compact subset of $\R^n$ is called a \textit{Kakeya set} is if it contains a line segment of length $1$ in every possible direction. To avoid ambiguity, we describe this more formally as follows.

\begin{definition} We say a compact set $K \subseteq \R^n$ is a \textit{Kakeya set} if for all directions $\omega \in S^{n-1}$, there exists a position $a \in \R^n$ such that the line segment $\ell_{\omega, a} := \{a + t \omega : 0 \leq t \leq 1\}$ satisfies $\ell_{\omega, a} \subseteq K$. 

\end{definition}

\begin{figure}
 \centering
 
  \begin{subfigure}[b]{0.3\textwidth}
          \centering
\includegraphics[width=0.9\linewidth]{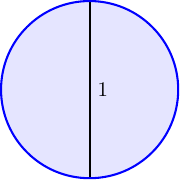}
          \caption{Area: $\frac{\pi}{4} = 0.785\dots$}
          
     \end{subfigure}
    \quad
\begin{subfigure}[b]{0.3\textwidth}   
\centering
 \includegraphics[width=0.9\linewidth]{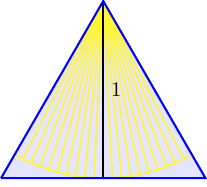}
      \caption{Area: $\frac{1}{\sqrt{3}} = 0.577\dots$}
          
     \end{subfigure} 
     \quad
\begin{subfigure}[b]{0.3\textwidth}
          \centering
\includegraphics[width=0.9\linewidth]{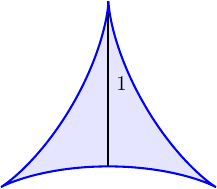}
          \caption{{Area: $\frac{\pi}{8} = 0.392\dots$}}
          \label{fig:Kak ex C}
     \end{subfigure}\\
			\caption{Three examples of Kakeya sets.}
\label{fig: Kak ex}
\end{figure}

In Figure~\ref{fig: Kak ex} we illustrate three simple examples of planar Kakeya sets: a disc, an equilateral triangle and a deltoid. The most instructive is the equilateral triangle. If we imagine pinning one end of a line segment to the top corner, we can then sweep it along a $60^{\circ}$ angle whilst staying within the triangle. This gives $60^{\circ}$ worth of directions pointing upwards and $60^{\circ}$ worth of directions pointing downwards, so a total of $120^{\circ}$ of directions. We can do this for all three corners of the triangle to account for all $360^{\circ}$.

If we compare the areas of the three Kakeya sets in Figure~\ref{fig: Kak ex}, we see they get successively smaller. In fact, the deltoid has half the area of the disc. This means that, within the deltoid, we have managed to pack the line segments much more tightly. In the early twentieth century, S\=oichi Kakeya was interested in these shapes and understanding different ways to pack lines into small spaces. From his work, a natural question arose: what is the smallest possible area of a planar Kakeya set?

\subsection{Besicovitch's surprising answer}\label{sec: Besicovitch construction} Around the same time, Besicovitch was independently thinking about similar questions. He was able to devise new and surprising ways to pack lines, which produced planar Kakeya sets with much, much smaller area than any example seen before. In fact, he was able to show that
\begin{center}
there exist planar Kakeya sets of \textit{arbitrarily small} area!     
\end{center}
That is, for any $\varepsilon > 0$, there exists some Kakeya set $K \subseteq \R^2$ with $|K| < \varepsilon$. Here and below, given $E \subseteq \R^d$ for $d \in \N$, we let $|E|$ denote the $d$-dimensional measure of $E$. 

\begin{figure}
 \centering
 
  \begin{subfigure}[b]{0.45\textwidth}
          \centering
\includegraphics[scale=1.2]{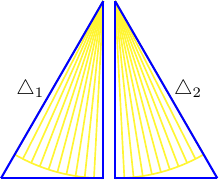}
          \caption{Cut $\triangle$ into subtriangles $\triangle_1$, $\triangle_2$.}
          
     \end{subfigure}
    \quad
\begin{subfigure}[b]{0.45\textwidth}   
\centering
 \includegraphics[scale=1.2]{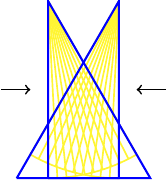}
          \caption{Horizontally shift $\triangle_1$, $\triangle_2$.}
          
     \end{subfigure}\\
			\caption{The cut-and-shift procedure.}
\label{fig: cut_and_shift}
\end{figure}

How is this possible? How can we pack so many line segments into such a tiny space?\footnote{Here we sketch a simplification of Besicovitch's argument. A detailed history of the problem, together with the precise description of the method, can be found in \cite[Chapter 7]{Falconer}.} To understand this, we go back to the example of the equilateral triangle, which we denote by $\triangle$. In Figure~\ref{fig: cut_and_shift}, we carry out the following `cut-and-shift' procedure to arrive at a new shape:
\begin{itemize}
    \item[(a)] Cut $\triangle$ into two `subtriangles' $\triangle_1$ and $\triangle_2$ along the vertical midline.
    \item[(b)] Shift $\triangle_1$ and $\triangle_2$ horizontally in opposite directions so that they overlap.
\end{itemize}
The key observation is that, as we cut-and-shift, we do not lose any of the line segments swept out from the top corner. Furthermore, because of the overlap between the shifted subtriangles, the new shape has smaller area than $\triangle$. The procedure therefore results in a shape with smaller area than the equilateral triangle $\triangle$, but which still contains line segments in $120^{\circ}$ of directions. This on its own is not very impressive. However, if we apply the procedure over and over many times, then we can amplify its effect and arrive at a startling conclusion.

\begin{figure}
 \centering
  \begin{subfigure}[htp]{0.46\textwidth}
          \centering
\includegraphics[scale=1.2]{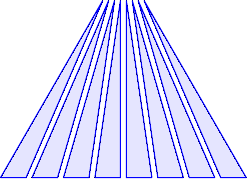}
          \caption{Cut $\triangle$ into $2^m$ subtriangles.}
          
     \end{subfigure}
    \quad
\begin{subfigure}[htp]{0.46\textwidth}   
\centering
 \includegraphics[scale=1.2]{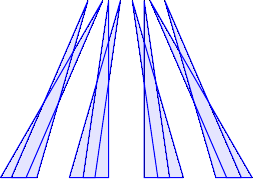}
          \caption{Horizontally shift to form $2^{m-1}$ pairs.}
     \end{subfigure}\\
     \medskip

  \begin{subfigure}[htp]{0.46\textwidth}
          \centering
\includegraphics[scale=1.2]{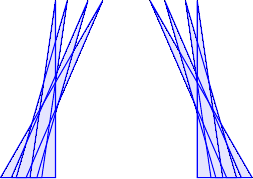}
        \caption{Form $2^{m-2}$ pairs of pairs...}
     \end{subfigure}
    \quad
\begin{subfigure}[htp]{0.46\textwidth}   
\centering
 \includegraphics[scale=1.2]{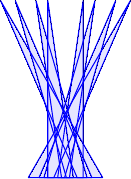}
          \caption{The final shape $K$.}
          \label{fig: final Perron}
     \end{subfigure}\\ 
     
			\caption{Perron tree construction of a Kakeya set.}
\label{fig: Perron}
\end{figure}

Rather than cutting $\triangle$ into just two subtriangles, in Figure~\ref{fig: Perron} we cut it into $2^m$ subtriangles for some large value of $m \in \N$. As illustrated, we then shift the subtriangles horizontally, first to form pairs, and then to form pairs of pairs, and then pairs of pairs of pairs and so on. By doing so, we ensure that there is a huge amount of overlap. At the end of this process, we arrive at a shape $K$ called a \textit{Perron tree}, represented in Figure~\ref{fig: final Perron}. The huge overlap means that the Perron tree $K$  has much, much smaller area than the original triangle $\triangle$. By choosing the number $2^m$ of subtriangles to be very large, and being careful how we shift, we can ensure that the area of $K$ is as small as we wish. As before, $K$ contains (shifted versions of) all the line segments spanned by the top corner of $\triangle$, so lines in $120^{\circ}$ of directions. We can then take three rotated copies of $K$ to form a true Kakeya set containing line segments pointing in all directions.\footnote{Taking three copies of $K$ will increase the area by a factor of $3$, but we can compensate for this by ensuring that the area of $K$ is three times smaller than our target value.}

The Perron tree construction produces a Kakeya set with area smaller than any prescribed value. In fact, combining these ideas with a limiting argument, Besicovitch was able to go one step further and show there exist Kakeya sets with \textit{zero} area!

\begin{theorem}[Besicovitch] There exists a Kakeya set $K \subseteq \R^2$ with $|K| = 0$.
\end{theorem}

Given a Kakeya set $K \subseteq \R^2$, the product $K \times [0,1]^{n-2}$ is easily seen to be a Kakeya set in $\R^n$. Hence, Besicovitch's theorem also implies that for every $n \geq 2$ there exists a Kakeya set $K \subseteq \R^n$ with $|K| = 0$.

Besicovitch's construction is undoubtedly striking and beautiful. However, around the time of their discovery, it is perhaps fair to say that measure zero Kakeya sets were something of a curiosity. Besicovitch did apply his construction to solve a problem in integration theory, but beyond this it did not appear to have much use. 

Mathematical ideas have a wonderful habit of reappearing in new and seemingly unrelated contexts. What was once thought a curiosity can later emerge as a profoundly important idea. Kakeya sets are a striking example of this. Fifty years after Besicovitch's work, startling discoveries uncovered their true significance. 

%%%%%%%%%%%%%%%%%%%%%%%%%%%%%%%%%%%%%%%%%%%%%%%%%%%%%%%%%%%%%%%%%%%%%%%%%%%%%%%%%%%%%%%%%%%%%%%%

%    Convergence of Fourier integrals

%%%%%%%%%%%%%%%%%%%%%%%%%%%%%%%%%%%%%%%%%%%%%%%%%%%%%%%%%%%%%%%%%%%%%%%%%%%%%%%%%%%%%%%%%%%%%%%%

\section{Convergence of Fourier integrals}

%%%%%%%%%%%%%%%%%%%%%%%%%%%%%%%%%%%%%%%%%%%%%%%%%%%%%%%%%%%%%%%%%%%%%%%%%%%%%%%%%%%%%%%%%%%%%%%%

%    Fourier inversion

%%%%%%%%%%%%%%%%%%%%%%%%%%%%%%%%%%%%%%%%%%%%%%%%%%%%%%%%%%%%%%%%%%%%%%%%%%%%%%%%%%%%%%%%%%%%%%%%

\subsection{Fourier inversion} Fourier analysis is concerned with representing functions in terms of sine and cosine waves (or, equivalently, complex exponentials). Ostensibly, this has little to do with the kind of elementary geometric problems discussed above. However, as we shall see, it is within the world of Fourier analysis that the true significance of Kakeya sets and Besicovitch's construction comes to light. 

We begin by recalling some of the basic principles of the Fourier theory. Given $f \in L^1(\R^n)$, the \textit{Fourier transform} is defined by\footnote{Recall that $L^p(\R^n)$ denotes the space of all measurable functions $f \colon \R^n \to \C$ satisfying $\|f\|_{L^p(\R^n)} := \big(\int_{\R^n} |f|^p \big)^{1/p} < \infty$.}
\begin{equation}\label{eq: Fourier transform}
    \widehat{f}(\xi) := \int_{\R^n} e^{-2\pi i x \cdot \xi} f(x)\,\ud x, \qquad \xi \in \widehat{\R}^n.
\end{equation}
Here we can think of the $\widehat{\R}^n$ notation in \eqref{eq: Fourier transform} as simply denoting a copy of $\R^n$. We include the hat to remind ourselves this is the domain of the Fourier transform, since later it will be useful to distinguish the domain of $\widehat{f}$ and the domain of $f$.

Roughly speaking, $\widehat{f}(\xi)$ corresponds to the amplitude of the component of $f$ which oscillates with fixed frequency vector $\xi \in \widehat{\R}^n$. This leads to the fundamental \textit{Fourier inversion formula}.

\begin{theorem}[Fourier inversion formula]\label{thm: inversion} If $f \in L^1(\R^n)$ \textbf{and} $\widehat{f} \in L^1(\widehat{\R}^n)$, then
\begin{equation}\label{eq: Fourier inversion}
    f(x) = \int_{\R^n} e^{2 \pi i x \cdot \xi} \widehat{f}(\xi)\,\ud \xi \qquad \textrm{for almost every $x \in \R^n$.}
\end{equation}
\end{theorem}

\begin{figure}
    \centering
    \includegraphics[width=0.4\linewidth]{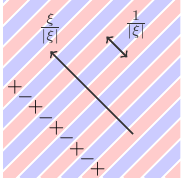}
    \caption{The level sets of the plane wave $x \mapsto \cos(2 \pi x \cdot \xi)$, corresponding to the real part of $e^{2 \pi i x \cdot \xi}$. The alternating bands represent regions where the function is either positive or negative. The wave has frequency $|\xi|$, which means that the wavelength (distance between consecutive peaks) is $1/|\xi|$. Each wave also has an orientation: the peaks and troughs lie orthogonal to $\xi$.}
    \label{fig: plane wave}
\end{figure}

The inversion formula is one realisation of the fundamental principle of Fourier analysis. In particular, under the above hypotheses, \eqref{eq: Fourier inversion} represents $f$ as a superposition of simple waves 
\begin{equation*}
    x \mapsto e^{2 \pi i x \cdot \xi} \widehat{f}(\xi).
\end{equation*}
Each simple wave oscillates with a fixed frequency $|\xi|$ and has fixed amplitude $\widehat{f}(\xi)$. For dimensions $n \geq 2$, each wave also has an orientation: the peaks and troughs lie orthogonal to $\xi$. See Figure~\ref{fig: plane wave}. 

%%%%%%%%%%%%%%%%%%%%%%%%%%%%%%%%%%%%%%%%%%%%%%%%%%%%%%%%%%%%%%%%%%%%%%%%%%%%%%%%%%%%%%%%%%%%%%%%

%    Shortcomings of the inversion formula

%%%%%%%%%%%%%%%%%%%%%%%%%%%%%%%%%%%%%%%%%%%%%%%%%%%%%%%%%%%%%%%%%%%%%%%%%%%%%%%%%%%%%%%%%%%%%%%%

\subsection{Shortcomings of the inversion formula} The hypotheses of Theorem~\ref{thm: inversion} require both $f \in L^1(\R^n)$ and $\widehat{f} \in L^1(\widehat{\R}^n)$. These conditions are natural: we need $f \in L^1(\R^n)$ to define the integral in \eqref{eq: Fourier transform} and we need $\widehat{f} \in L^1(\widehat{\R}^n)$ to define the integral in \eqref{eq: Fourier inversion}. Despite this, the condition $\widehat{f} \in L^1(\widehat{\R}^n)$ is, in some sense, restrictive and a severe limitation on the theory. It requires
\begin{equation*}
    \int_{\widehat{\R}^n} |\widehat{f}(\xi)| \,\ud \xi < \infty,
\end{equation*}
which means $\widehat{f}(\xi)$ has fast decay as $|\xi| \to \infty$. In other words, the high frequency components of $f$ must have very small amplitude. Intuitively, the graph of such a function cannot have sudden jumps, since jumps require high frequencies. This forces $f$ to be rather regular, which drastically limits the applicability of Theorem~\ref{thm: inversion}. Moreover, often in applications of Fourier analysis, it is critical to understand precisely \textit{high} frequency behaviour. Since the inversion formula does not apply to functions with significant high frequency components, it is poorly suited to many of these applications.

 %%%%%%%%%%%%%%%%%%%%%%%%%%%%%%%%%%%%%%%%%%%%%%%%%%%%%%%%%%%%%%%%%%%%%%%%%%%%%%%%%%%%%%%%%%%%%%%%

%    Convergence of Fourier integrals: one dimensional case

%%%%%%%%%%%%%%%%%%%%%%%%%%%%%%%%%%%%%%%%%%%%%%%%%%%%%%%%%%%%%%%%%%%%%%%%%%%%%%%%%%%%%%%%%%%%%%%%

\subsection{Convergence of Fourier integrals: one dimensional case}  Without the strong hypothesis $\widehat{f} \in L^1(\widehat{\R}^n)$, the integral formula in \eqref{eq: Fourier inversion} is not well defined. Thus, on the face of it, it is unclear how to make sense of Fourier inversion without requiring $\widehat{f} \in L^1(\widehat{\R}^n)$. However, restricting to the $n=1$ case, we can still make sense of the \textit{partial Fourier integrals}
\begin{equation}\label{eq: S_R}
    S_Rf(x) := \int_{-R}^R e^{2 \pi i x \xi} \widehat{f}(\xi)\,\ud \xi, \qquad R > 0. 
\end{equation}
Since the integration now takes place over a bounded interval, the integral in \eqref{eq: S_R} always converges. Moreover, with a small amount of work, \eqref{eq: S_R} can be rewritten as the convolution
\begin{equation}\label{eq: conv formula}
    S_Rf(x) = D_R \ast f(x) \qquad \textrm{where} \qquad D_R(x) := \frac{\sin (2\pi R x)}{\pi x}.
\end{equation}
The precise formula for $D_R$ is not too important for us. However, what is important is that \eqref{eq: conv formula} makes no reference to the Fourier transform of $f$ and, consequently, one can easily show $S_Rf$ is well defined for any $f \in L^p(\R)$ for $1 \leq p < \infty$ (in particular, we no longer require $f \in L^1(\R)$). It is natural to ask whether, at this level of generality, $f$ can be reconstructed from the partial Fourier integrals $S_Rf$. 

\begin{question}\label{qtn: conv} Does $S_Rf \to f$ as $R \to \infty$?
\end{question}

A positive answer to Question~\ref{qtn: conv} acts as a surrogate for the inversion formula \eqref{eq: Fourier inversion}. However, as stated, the question is somewhat vague. Since we are dealing with limits of families of \textit{functions}, there are many different modes of convergence to consider. The following theorem, which combines two foundational results in Fourier analysis, treats the two most natural modes. 

\begin{theorem}[Riesz \cite{Riesz1928}, Carleson--Hunt \cite{Carleson1966, Hunt1968}]\label{thm: 1d Fourier} Let $1 < p < \infty$ and $f \in L^p(\R)$. Then
\begin{enumerate}[i)]
    \item $\|S_Rf - f\|_{L^p(\R)} \to 0$ as $R \to \infty$;
    \item $S_Rf(x) \to f(x)$ as $R \to \infty$ for almost every $x \in \R$.
\end{enumerate}
\end{theorem}

Theorem~\ref{thm: 1d Fourier} represents a sweeping and powerful extension of the inversion formula in the one-dimensional setting. Since the hypotheses only require $f \in L^p(\R)$ (and no additional conditions on the Fourier transform), the results apply in great generality, including to very rough functions. The first part of Theorem~\ref{thm: 1d Fourier}, due to Riesz~\cite{Riesz1928}, shows that $S_Rf$ converges to $f$ with respect to the metric structure on $L^p(\R)$. Since the $L^p$-norms involve an integral, this can be interpreted as showing convergence holds `on average'.  The second result, on almost everywhere convergence, involves no averaging and is far more subtle. This is the celebrated Carleson--Hunt theorem \cite{Carleson1966, Hunt1968}, a landmark in twentieth century analysis.

%%%%%%%%%%%%%%%%%%%%%%%%%%%%%%%%%%%%%%%%%%%%%%%%%%%%%%%%%%%%%%%%%%%%%%%%%%%%%%%%%%%%%%%%%%%%%%%%

%    Convergence of Fourier integrals: higher dimensions

%%%%%%%%%%%%%%%%%%%%%%%%%%%%%%%%%%%%%%%%%%%%%%%%%%%%%%%%%%%%%%%%%%%%%%%%%%%%%%%%%%%%%%%%%%%%%%%%

\subsection{Convergence of Fourier integrals: higher dimensions} 

As soon as we move to higher dimensions, convergence of Fourier integrals becomes a much more difficult and subtle problem. To begin with, there is no longer a `canonical' choice of truncation for defining the partial Fourier integrals $S_R$ as in \eqref{eq: S_R}. For instance, two natural choices are given by the \textit{square} and \textit{ball} partial integrals\footnote{Here $B(0,R)$ denotes the Euclidean ball centred at the origin of radius $R$.}
\begin{equation*}
    S_R^{\mathrm{square}}f(x) := \int_{[-R, R]^n} e^{2 \pi i x \cdot \xi} \widehat{f}(\xi)\,\ud \xi \quad \textrm{and} \quad S_R^{\mathrm{ball}}f(x) := \int_{B(0,R)} e^{2 \pi i x \cdot \xi} \widehat{f}(\xi)\,\ud \xi,
\end{equation*}
respectively. The above definitions coincide for $n = 1$, but are different beasts in higher dimensions. Taking the limit of the partial Fourier integrals $S_R^{\mathrm{square}}f$ and $S_R^{\mathrm{ball}}f$ as $R \to \infty$ corresponds to two different \textit{summation methods} for $f$: that is, different ways to attempt to reconstruct $f$ from its Fourier transform.

\medskip \noindent \textit{Square summation.} The summation method for the square partial integrals is relatively easy to understand. In this setting, the one-dimensional theory from Theorem~\ref{thm: 1d Fourier} extends directly to higher dimensions.

\begin{theorem}\label{thm: sq Fourier} For $n \geq 1$, let $1 < p < \infty$ and $f \in L^p(\R^n)$. Then
\begin{enumerate}[i)]
    \item $\|S_R^{\mathrm{square}}f - f\|_{L^p(\R^n)} \to 0$ as $R \to \infty$;
    \item $S_R^{\mathrm{square}}f(x) \to f(x)$ as $R \to \infty$ for almost every $x \in \R^n$.
\end{enumerate}
\end{theorem}

Since the square $[-R, R]^n$ is formed by a product of intervals $[-R,R]$, the square partial integral operators $S_R^{\mathrm{square}}$ have an underlying tensor product structure. Using this, the norm convergence problem for $S_R^{\mathrm{square}}$ immediately reduces to the $n=1$ case. A similar, but this time somewhat non-trivial, reduction applies to almost everywhere convergence \cite{Fefferman1971a}. In this way, Theorem~\ref{thm: sq Fourier} can be derived as a corollary of Theorem~\ref{thm: 1d Fourier}.

\medskip \noindent \textit{Ball summation.} In contrast with the previous case, the ball partial integral operators $S_R^{\mathrm{ball}}$ do not possess a tensor structure. Consequently, there is no obvious or easy way to reduce  the convergence problem for the ball summation method to the $n=1$ case. This situation is `genuinely' higher dimensional. 

Understanding the norm convergence problem for $S_R^{\mathrm{ball}}$ was a major open problem in the 1960s. Using Plancherel's theorem, it is a trivial matter to see that norm convergence holds for $p = 2$; that is, 
\begin{equation*}
    \|S_R^{\mathrm{ball}} f - f\|_{L^2(\R^n)} \to 0 \qquad \textrm{as $R \to \infty$} \qquad \textrm{for all $f \in L^2(\R^n)$. }
\end{equation*}
However, the behaviour for $p \neq 2$ is much less clear. It came as a great surprise when, in 1971, Charles Fefferman \cite{Fefferman1971} produced a counterexample showing that norm convergence holds \textit{only} in the trivial $p=2$ case!

\begin{theorem}[Fefferman \cite{Fefferman1971}]\label{thm: Fefferman} For $n \geq 2$ and $1 \leq p \leq \infty$ with $p \neq 2$, there exists some $f \in L^p(\R^n)$ such that 
\begin{equation*}
    \|S_R^{\mathrm{ball}}f - f\|_{L^p(\R^n)} \not\to 0 \qquad \textrm{as $R \to \infty$.} 
\end{equation*}  
In particular, norm convergence holds \textbf{only} for $p = 2$. 
\end{theorem}

Fefferman's theorem also rules out almost everywhere convergence for $p \neq 2$.
%\footnote{However, for $p = 2$, the almost everywhere convergence problem for $S_R^{\mathrm{ball}}$ remains an open question, even in the plane! In particular, it is not known whether $S^{\mathrm{ball}}_Rf(x) \to f(x)$ as $R \to \infty$ for almost every $x \in \R^2$ whenever $f \in L^2(\R^2)$. This problem seems well beyond the reach of current technology.} 
These results illustrate a remarkable difference between the square and ball summation methods. Whilst the square Fourier integrals have a robust convergence theory and behave just about as well as one could hope for, the ball Fourier integrals behave just about as badly as one could dread!

Fefferman's theorem is celebrated not only for its startling conclusion, but also for the powerful ideas used in the proof. In particular, his counterexample relies heavily on Besicovitch's construction of a measure zero Kakeya set! Moreover, as we shall see, the proof of Theorem~\ref{thm: Fefferman} revealed a profound relationship between the geometry of Kakeya sets and the theory of the Fourier transform in two and higher dimensions. This marked a tremendous shift in perspective, and the ramifications of Fefferman's insight remain a central and intensely active topic of contemporary research in the modern Fourier analysis and beyond.

%%%%%%%%%%%%%%%%%%%%%%%%%%%%%%%%%%%%%%%%%%%%%%%%%%%%%%%%%%%%%%%%%%%%%%%%%%%%%%%%%%%%%%%%%%%%%%%%

%    Fefferman's counterexample: why Kakeya sets are important
%%%%%%%%%%%%%%%%%%%%%%%%%%%%%%%%%%%%%%%%%%%%%%%%%%%%%%%%%%%%%%%%%%%%%%%%%%%%%%%%%%%%%%%%%%%%%%%%

\section{Fefferman's counterexample: why Kakeya sets are important}\label{sec: Fefferman}

Here we provide a non-rigorous sketch of the proof of Theorem~\ref{thm: Fefferman}, focusing on how Kakeya sets arise in Fefferman's counterexample. We shall focus on the planar case $n=2$, which already contains all the essential ideas. Moreover, once the $n = 2$ case is known, standard arguments easily lift the result to higher dimensions.

%%%%%%%%%%%%%%%%%%%%%%%%%%%%%%%%%%%%%%%%%%%%%%%%%%%%%%%%%%%%%%%%%%%%%%%%%%%%%%%%%%%%%%%%%%%%%%%%

%    A low pass filter

%%%%%%%%%%%%%%%%%%%%%%%%%%%%%%%%%%%%%%%%%%%%%%%%%%%%%%%%%%%%%%%%%%%%%%%%%%%%%%%%%%%%%%%%%%%%%%%%

\subsection{A low pass filter} Consider the partial integral operator $S := S_1^{\mathrm{ball}}$, given by
\begin{equation*}
    Sf(x) = \int_{B(0,1)} e^{2 \pi i x \cdot \xi} \widehat{f}(\xi)\,\ud \xi. 
\end{equation*}
This has a natural interpretation as a \textit{low pass filter}. It acts on an input function $f$ by removing all high frequency components, which corresponds to multiplying $\widehat{f}$ against the characteristic function $\chi_{B(0,1)}$ of the unit ball. A basic question is to understand the continuity of this filtering operation. Given $1 \leq p \leq \infty$, the continuity of $S$ on the metric space $L^p(\R^2)$ is equivalent to the existence of a constant $C > 0$ such that
\begin{equation}\label{eq: ball multiplier}
    \|Sf\|_{L^p(\R^2)} \leq C\|f\|_{L^p(\R^2)} \qquad \textrm{for all $f \in L^p(\R^2)$.}
\end{equation}

The question of whether $S$ is continuous on $L^p(\R^2)$ is closely related to the question of whether $L^p$ convergence
\begin{equation}\label{eq: ball Lp conv}
    \|S_R^{\mathrm{ball}}f - f\|_{L^p(\R^2)} \to 0 \qquad \textrm{as $R \to \infty$} \qquad \textrm{for all $f \in L^p(\R^2)$}
\end{equation}
holds. To understand this relationship, note:
\begin{enumerate}[i)]
    \item $L^p$ convergence always holds for smooth, rapidly decaying functions $g$. Indeed, in this case $g$ and $\widehat{g}$ are both well behaved and one may establish convergence via elementary tools such as the inversion formula (Theorem~\ref{thm: inversion}).
    \item Any $f \in L^p(\R^2)$ can be approximated to arbitrary accuracy in $L^p(\R^2)$ by some smooth, rapidly decaying function $g \in L^p(\R^2)$. 
\end{enumerate}
Together, i) and ii) tell us that the limit identity in \eqref{eq: ball Lp conv} holds for all $f$ belonging to a dense subspace of $L^p(\R^2)$. If the operators $(S_R^{\mathrm{ball}})_{R > 0}$ were (equi)continuous, then one could take limits and pass from the dense subspace to \textit{all} $L^p(\R^2)$ functions. More precisely, suppose we have an inequality of the form
\begin{equation}\label{eq: equicontinuity}
    \|S_R^{\mathrm{ball}}f\|_{L^p(\R^2)} \leq C\|f\|_{L^p(\R^2)} \qquad \textrm{for all $f \in L^p(\R^2)$,}
\end{equation}
where the constant $C > 0$ is uniform in $R$. This is a statement about the continuity of the family of linear operators $(S_R^{\mathrm{ball}})_{R > 0}$ on $L^p(\R^2)$. Following the above discussion, the inequality \eqref{eq: equicontinuity} implies the $L^p$ convergence result \eqref{eq: ball Lp conv}. 

It turns out that the converse of the above observation is also true. That is, if $L^p$ convergence result \eqref{eq: ball Lp conv} holds, then the inequality \eqref{eq: equicontinuity} must also be true. This implication is deeper, relying on the \textit{principle of uniform boundedness}; we do not discuss the details here. The key takeaway is that, in order to prove Fefferman's theorem that the $L^p$ convergence result \eqref{eq: ball Lp conv} fails, it suffices to show that the inequality \eqref{eq: equicontinuity} fails for all $p \neq 2$. We shall demonstrate this failure for $R = 1$, in which case \eqref{eq: equicontinuity} agrees with \eqref{eq: ball multiplier}. Thus, matters are reduced to showing that 
\begin{equation*}
    \textrm{the inequality \eqref{eq: ball multiplier} fails for all $p \neq 2$. }
\end{equation*}
This is the goal for the remainder of the section. 
%%%%%%%%%%%%%%%%%%%%%%%%%%%%%%%%%%%%%%%%%%%%%%%%%%%%%%%%%%%%%%%%%%%%%%%%%%%%%%%%%%%%%%%%%%%%%%%%

%    Elementary properties of the Fourier transform
%%%%%%%%%%%%%%%%%%%%%%%%%%%%%%%%%%%%%%%%%%%%%%%%%%%%%%%%%%%%%%%%%%%%%%%%%%%%%%%%%%%%%%%%%%%%%%%%
\subsection{Elementary properties of the Fourier transform} Fefferman's argument relies on two elementary properties of the Fourier transform. 
%%%%%%%%%%%%%%%%%%%%%%%%%%%%%%%%%%%%%%%%%%%%%%%%%%%%%%%%%%%%%%%%%%%%%%%%%%%%%%%%%%%%%%%%%%%%%%%%

%    The uncertainty principle
%%%%%%%%%%%%%%%%%%%%%%%%%%%%%%%%%%%%%%%%%%%%%%%%%%%%%%%%%%%%%%%%%%%%%%%%%%%%%%%%%%%%%%%%%%%%%%%%

\begin{figure}
 \centering
\includegraphics[width=0.65\linewidth]{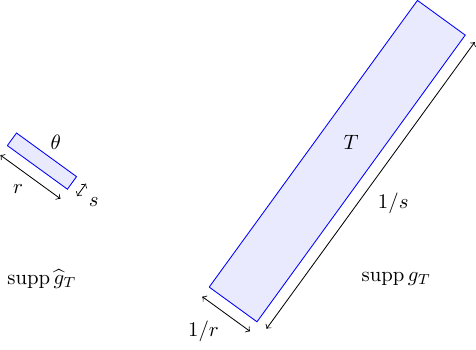}

\caption{The uncertainty principle in $\R^2$. The Fourier transform $\widehat{g}_T$ is a bump adapted to the $r \times s$ rectangle $\theta$. The function $g_T$ is a normalised bump adapted to the $1/r \times 1/s$ dual rectangle $T$. }
\label{fig: uncertainty}
\end{figure}

\medskip \noindent \textit{The uncertainty principle.} Rather than a single result or set of results, the uncertainty principle is a guiding heuristic which permeates all aspects of Fourier analysis. Loosely speaking, it says that:
\begin{quote}
    Localising $\widehat{f}$ in the frequency domain forces $f$ to spread out in the spatial domain.
\end{quote}
To illustrate this idea on the real line, for a small number $r > 0$, consider the scaled Gaussian $f_r(x) := r e^{- \pi r^2 x^2}$ with variance $1/r^2$. Then $\widehat{f}_r(\xi) = e^{-\pi x^2/r^2}$ is a Gaussian with variance $r^2$. The variance of $f_r$ and the variance $\widehat{f}_r$ satisfy a reciprocal relationship. By shrinking $r$, the distribution of $\widehat{f}_r$ becomes narrower (more certain) whilst the distribution of $f_r$ becomes wider (more uncertain). This behaviour is a manifestation of the uncertainty principle. 

We are interested in a two-dimensional variant of the above example. For small numbers $r$, $s > 0$, let $\theta \subseteq \R^2$ be a rectangle of dimensions $r \times s$, with arbitrary centre and orientation. We define the \textit{dual rectangle} $T \subseteq \R^2$ to be the $1/r \times 1/s$ rectangle centred at the origin with the same orientation as $\theta$. We illustrate this definition in Figure~\ref{fig: uncertainty}.

Suppose $g_T$ is a function on $\R^2$ whose Fourier transform $\widehat{g}_T$ is a \textit{bump adapted to the rectangle} $\theta$. By this, we  mean that $\widehat{g}_T$ looks like\footnote{We intentionally keep the discussion vague here to avoid technicalities. For the purposes of this exposition, one could simply think of $\widehat{g}_T = \chi_{\theta}$. However, in practice, it is useful to take a smoothed out version of the characteristic function, since the Fourier transform is poorly behaved in the presence of discontinuities.} a smoothed out version of the characteristic function $\chi_{\theta}$. Then the uncertainty principle tells us that the support of $g_T$ spreads out over the dual rectangle $T$. This is a two-dimensional analogue of the Gaussian example above. In particular, if we think of $\widehat{g}_T$ as a smoothed out version of the characteristic function $\chi_{\theta}$ of $\theta$, then, at least morally, $|g_T|$ should be a smoothed out version of the normalised characteristic function $|T|^{-1} \chi_T$.\footnote{The normalising factor $|T|^{-1}$ naturally appears owing to the scaling properties of the Fourier transform. It plays the same role as the $r$ factor in the definition of $f_r$ in the Gaussian example.}\footnote{This discussion applies to the \textit{absolute value} $|g_T|$. The function $g_T$ itself carries some oscillation, which is dictated by the location of $\theta$ in the frequency domain.}

%%%%%%%%%%%%%%%%%%%%%%%%%%%%%%%%%%%%%%%%%%%%%%%%%%%%%%%%%%%%%%%%%%%%%%%%%%%%%%%%%%%%%%%%%%%%%%%%

%    Modulation vs translation symmetry

%%%%%%%%%%%%%%%%%%%%%%%%%%%%%%%%%%%%%%%%%%%%%%%%%%%%%%%%%%%%%%%%%%%%%%%%%%%%%%%%%%%%%%%%%%%%%%%%
    
\medskip \noindent \textit{Translation vs modulation symmetry.} Given $y \in \R^2$, the Fourier transform intertwines the \textit{translation} and \textit{modulation} operations
\begin{equation*}
 (\tau_yf)(x) := f(x - y) \qquad \textrm{and} \qquad  \mathrm{Mod}_y h(\xi) := e^{-2\pi y \cdot \xi} h(\xi).
\end{equation*}
In particular, the integral formula \eqref{eq: Fourier transform}, together with a change of variables, implies 
\begin{equation*}
    (\tau_yf)\;\widehat{}\;(\xi) =  \int_{\R^n} e^{-2 \pi i x \cdot \xi} f(x-y) \,\ud x = e^{-2\pi i y \cdot \xi}\int_{\R^n} e^{-2 \pi i x \cdot \xi} f(x) \,\ud x = \mathrm{Mod}_y \widehat{f}(\xi).
\end{equation*}
Thus, translating the function $f$ in the physical domain corresponds to modulating $\widehat{f}$ in the frequency domain. 

%%%%%%%%%%%%%%%%%%%%%%%%%%%%%%%%%%%%%%%%%%%%%%%%%%%%%%%%%%%%%%%%%%%%%%%%%%%%%%%%%%%%%%%%%%%%%%%%

%    The counterexample

%%%%%%%%%%%%%%%%%%%%%%%%%%%%%%%%%%%%%%%%%%%%%%%%%%%%%%%%%%%%%%%%%%%%%%%%%%%%%%%%%%%%%%%%%%%%%%%%
    
\subsection{The counterexample} We now have all the ingredients we need to sketch Fefferman's argument. We break the proof into steps. 

\medskip \noindent \textit{Step 0: Recap.} We begin by recalling what it is we wish to prove. As above, $S$ denotes the low pass filter
\begin{equation*}
    Sf(x) := \int_{B(0,1)} e^{2 \pi i x \cdot \xi} \widehat{f}(\xi)\,\ud \xi.
\end{equation*}
This removes components of $f$ oscillating with frequency $|\xi| > 1$. By our earlier reductions, our goal is to show
\begin{equation}\label{eq: ball multiplier recall}
    \|Sf\|_{L^p(\R^2)} \leq C\|f\|_{L^p(\R^2)} \qquad \textrm{for all $f \in L^p(\R^2)$}
\end{equation}
\textbf{fails} whenever $p \neq 2$. Suppose, to the contrary, that the inequality \eqref{eq: ball multiplier recall} holds. What does this tell us? 

The $L^p$ norms are sensitive to the distributional properties (that is, the shape of the graph) of a given function. Since $L^{\infty}$ functions are by definition bounded, their distributions are typically broad and flat. On the other hand, $L^1$ functions are defined in terms of the convergence of an integral, which requires some decay at $\infty$, so their distributions are typically tall and narrow. A function in $L^p$ for some $1 < p < \infty$ lies in a spectrum between these two extreme cases. 

If the inequality \eqref{eq: ball multiplier recall} were to hold, say for large $p > 2$, then this would tell us that the filter $S$ preserves certain distributional properties of the input. In particular, if the input function $f$ is broad and flat, then the output function $Sf$ should also be broad and flat. Therefore, to disprove the inequality, we shall construct an input function $f$ whose distributional properties change significantly under filtering:
\begin{itemize}
    \item The input $f$ will have a broad, flat distribution;
    \item The output $Sf$ will have a tall, narrow distribution. 
\end{itemize}
This kind of behaviour is incompatible with \eqref{eq: ball multiplier recall} for $p > 2$ and therefore provides a counterexample for this exponent regime. A similar dual construction can be used to rule out \eqref{eq: ball multiplier recall} for $p < 2$. 

\begin{figure}
 \centering
 
  \begin{subfigure}[m]{0.45\textwidth}
          \centering
\includegraphics[width=0.9\linewidth]{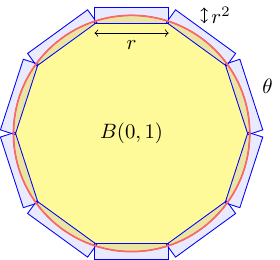}
          \caption{Schematic for $\widehat{g} = \sum_T \widehat{g}_T$. Each $\widehat{g}_T$ is supported on some $r \times r^2$ rectangle $\theta$ tangent to the boundary of $B(0,1)$.}
    \label{fig: Fefferman 1 a}      
     \end{subfigure}
    \quad
\begin{subfigure}[m]{0.45\textwidth}   
\centering
 \includegraphics[width=0.9\linewidth]{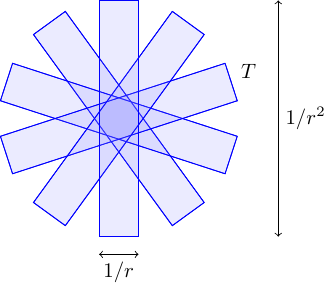}
          \caption{Schematic for $g = \sum_T g_T$. Each $g_T$ is essentially supported on a $1/r \times 1/r^2$ rectangle $T$ through the origin.}
     \label{fig: Fefferman 1 b}        
     \end{subfigure} 
   
\caption{The function $g$ from Step 2 of Fefferman's construction, viewed from both the frequency and physical domain. Each spatial rectangle $T$ in Figure~\ref{fig: Fefferman 1 b} is dual to some frequency rectangle $\theta$ in Figure~\ref{fig: Fefferman 1 a}.}
\label{fig: Fefferman 1}
\end{figure}

\medskip \noindent \textit{Step 1: Arrange frequency support.} We now turn to constructing the counterexample. That is, we set about finding a function $f$ with the above properties. Fix a small parameter $0 < r < 1$. As illustrated in Figure~\ref{fig: Fefferman 1 a}, we place $r \times r^2$ rectangles $\theta$ tangentially to the unit circle $\partial B(0,1)$ in the frequency domain $\widehat{\R}^2$. The choice of side lengths $r \times r^2$ reflects the curvature of the circle: with this choice, the $\theta$ can be arranged to fit neatly around the circle and form an efficient covering. 

\medskip \noindent \textit{Step 2: Initial physical support.} Given $\theta$, let $T$ denote the corresponding $1/r \times 1/r^2$ dual rectangle, which is centred at the origin. The collection of all such dual rectangles is illustrated in Figure~\ref{fig: Fefferman 1 b}. By \textbf{the uncertainty principle}, we can find a function $g_T$ such that $\widehat{g}_T$ is a bump adapted to $\theta$ and $g_T$ is a normalised bump adapted to $T$. We initially consider the function $g := \sum_T g_T$. The two sides of Figure~\ref{fig: Fefferman 1} then provide a schematic of the \textit{frequency support} of $g$ (the support of $\widehat{g}$ in $\widehat{\R}^2$) and \textit{physical support} of $g$ (the support of $g$ itself in $\R^2$).

\begin{figure}
 \centering
 
  \begin{subfigure}[b]{0.45\textwidth}
          \centering
\includegraphics[width=0.9\linewidth]{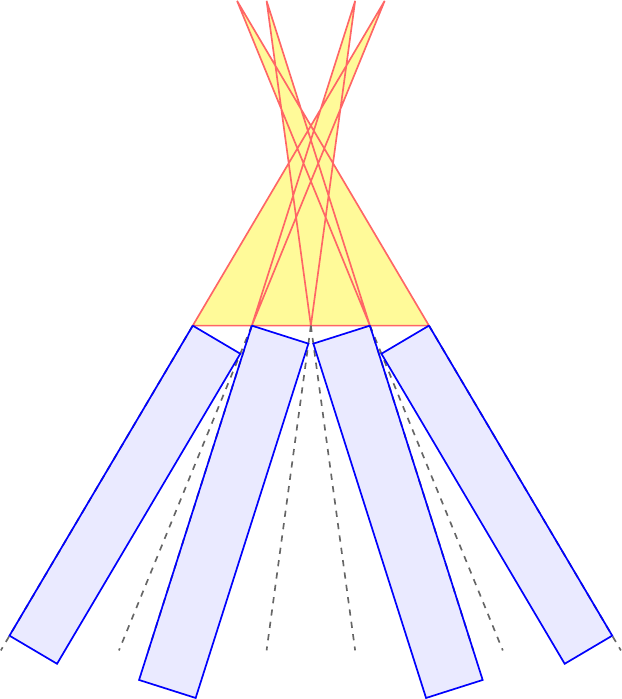}
          \caption{Schematic for $f = \sum_T f_T$, viewed in the physical domain. Each $f_T$ is essentially supported in a rectangle $T$ lying below the Perron tree. These rectangles are disjoint.}
    \label{fig: Fefferman 3 a}      
     \end{subfigure}
    \quad
\begin{subfigure}[b]{0.45\textwidth}   
\centering
 \includegraphics[width=0.9\linewidth]{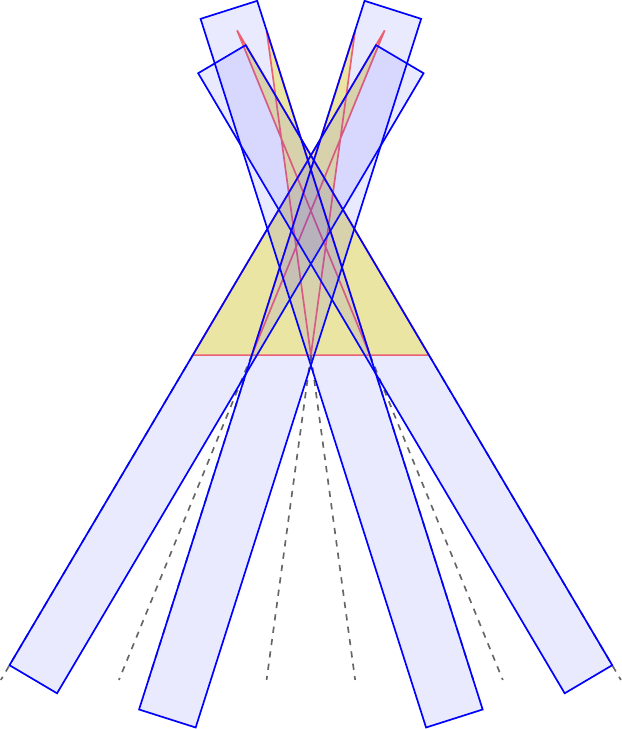}
          \caption{Schematic for $Sf = \sum_T Sf_T$. Morally, each $Sf_T$ is essentially supported on a double of $T$. The enlarged rectangles have huge overlap owing to the properties of the Perron tree. }
     \label{fig: Fefferman 3 b}        
     \end{subfigure} 
   
\caption{Applying the operator $S$ to the function $f$. The rectangles on the left represent the support of the input $f$. The rectangles on the right represent the support of the output $Sf$. }
\label{fig: Fefferman 3}
\end{figure}

\medskip \noindent \textit{Step 3: Rearrange physical support.} Using \textbf{translation vs modulation symmetry}, by independently modulating each $\widehat{g}_T$, we translate the $g_T$ in physical space. This corresponds to translating the underlying supporting rectangles $T$. We rearrange the $T$ so, rather than all passing through the origin, they are disjoint. Moreover, as illustrated in Figure~\ref{fig: Fefferman 3 a}, we strategically place the $T$ below the base of a Perron tree.\footnote{Since a single Perron tree only accounts for lines in $120^{\circ}$ of directions, in practice we use multiple rotated Perron trees. However, for the sake of simplicity we ignore this minor detail.} We let $f_T$ denote the corresponding translate of $g_T$ and $f = \sum_T f_T$. Figure~\ref{fig: Fefferman 3 a} then provides a schematic of the physical support of $f$.

Recall that $\widehat{f}_T$ is obtained by modulating $\widehat{g}_T$. Modulation does not alter the support of a function so $\widehat{f}_T$ and $\widehat{g}_T$ share the same same support. Thus, Figure~\ref{fig: Fefferman 1 a} provides a schematic of the frequency support of both $g = \sum_T g_T$ and $f = \sum_T f_T$. 

\medskip \noindent \textit{Step 4: Apply the operator.} We now consider what happens when we apply the filter $S$ to the function $f$ constructed above. By linearity, $Sf = \sum_T Sf_T$ and so it suffices to study the individual $f_T$.

We first describe the action of $S$ from the perspective of the frequency domain.  Recall that $\widehat{f}_T$ is supported on a frequency rectangle $\theta$, arranged along the boundary of the ball as in Figure~\ref{fig: Fefferman 1 a}. The rectangle $\theta$ is tangential to the unit circle, and the filter $S$ acts by cutting off frequencies to the unit ball. Consequently, $S$ effectively cuts the support $\theta$ of $\widehat{f}_T$ in half. 

We now consider the action of $S$ from the perspective of the spatial domain. As discussed above, $S$ effectively halves the frequency support of $f_T$. By \textbf{the uncertainty principle}, halving the frequency support causes the physical support to double. Consequently, in the spatial domain, each $Sf_T$ spreads out over a double of $T$. This phenomenon is illustrated in Figure~\ref{fig: Fefferman 3}.

\medskip \noindent \textit{Step 5: Conclusion.} As we can see in Figure~\ref{fig: Fefferman 3 b}, the strategic placement of the rectangles $T$ causes the doubles to pile up on top of one another. Moreover, since the Perron tree packs lines in different directions into a tiny space, the pile up between the doubles is enormous. This pile up is key.\medskip

From the previous steps we know:
\begin{itemize}
    \item The disjoint rectangles in Figure~\ref{fig: Fefferman 3 a} correspond to the supports of the functions $f_T$ forming $f = \sum_T f_T$;
    \item The piled-up doubles in Figure~\ref{fig: Fefferman 3 b} correspond to the supports of the functions $Sf_T$ forming $Sf = \sum_T Sf_T$. 
\end{itemize}
 Comparing both sides of Figure~\ref{fig: Fefferman 3}, we pass from an input function $f$ with a broad, flat distribution to an output function $Sf$ with a tall, narrow distribution. Thus, $S$ fails to preserve the distributional properties of the input, and therefore cannot be bounded on $L^p$.\medskip

The above argument is a rough sketch, and for the sake of brevity omits some details.\footnote{A significant omission is the effect of the oscillation of the individual terms of the sum $Sf = \sum_T Sf_T$. If there were a lot of destructive interference between the $Sf_T$, then this could, in principle, cancel out the effect of the pile up of the supports. However, it is reasonably straightforward to understand and control the oscillation in the proof.} Nevertheless, it serves to highlight the essential role of the geometry of the Kakeya sets in Fefferman's argument. It is only by arranging the rectangles in relation to (an approximation of) a measure zero Kakeya set that the pile up effect is sufficient to prove the theorem.

 Fefferman's profound discovery showed that the geometry of Kakeya sets underpins Fourier analysis in two and higher dimensions. Since Fourier analysis is so central to the study of many physical processes, it soon became apparent that Kakeya sets were also deeply connected to the theory of many important equations governing the laws of nature. This includes the wave equation, which describes a diverse range of phenomena such as water waves, seismic waves and sound waves, and also Schr\"odinger's famous equation from quantum mechanics \cite{Tao2001, Wolff1999}. Kakeya sets, which were once was a curiosity, suddenly became one of the most important objects in mathematical analysis! 

%%%%%%%%%%%%%%%%%%%%%%%%%%%%%%%%%%%%%%%%%%%%%%%%%%%%%%%%%%%%%%%%%%%%%%%%%%%%%%%%%%%%%%%%%%%%%%%%

%    Next steps in Fourier analysis

%%%%%%%%%%%%%%%%%%%%%%%%%%%%%%%%%%%%%%%%%%%%%%%%%%%%%%%%%%%%%%%%%%%%%%%%%%%%%%%%%%%%%%%%%%%%%%%%

\section{The next steps in Fourier analysis}

%%%%%%%%%%%%%%%%%%%%%%%%%%%%%%%%%%%%%%%%%%%%%%%%%%%%%%%%%%%%%%%%%%%%%%%%%%%%%%%%%%%%%%%%%%%%%%%%

%    A possible fix: the Bochner--Riesz Conjecture

%%%%%%%%%%%%%%%%%%%%%%%%%%%%%%%%%%%%%%%%%%%%%%%%%%%%%%%%%%%%%%%%%%%%%%%%%%%%%%%%%%%%%%%%%%%%%%%%

\subsection{A possible fix: the Bochner--Riesz Conjecture} 
\begin{figure}
    \centering
    \includegraphics[width=0.8\linewidth]{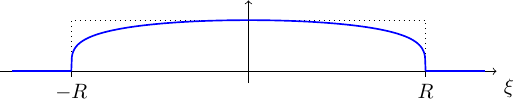}
    \caption{A radial slice of $m^{\alpha}_R(\xi) := (1 - |\xi|^2/R^2)^{\alpha}\chi_{B(0,R)}(\xi)$ from the definition of the Bochner--Riesz means. Here $\alpha = 1/4$. In general, for small $\alpha > 0$, the function is continuous but fails to be smooth around the boundary of the $B(0,R)$. }
    \label{fig: Bochner Riesz}
\end{figure}

Fefferman's example demonstrated a dramatic and highly unexpected failure of the ball Fourier summation method. In response to this, the hunt was on for a possible fix. A natural candidate is to consider slightly smoothed out versions of the ball summation method. Given $\alpha \geq 0$, we define the \textit{Bochner--Riesz means}
\begin{equation*}
    B_R^{\alpha}f(x) := \int_{B(0,R)} e^{2 \pi i x \cdot \xi} (1 - |\xi|^2/R^2)^{\alpha} \widehat{f}(\xi)\,\ud \xi, \qquad R > 0.
\end{equation*}
When $\alpha = 0$, this definition agrees with $S_R^{\mathrm{ball}}f$. However, for $\alpha > 0$, the Bochner--Riesz mean $B_R^{\alpha}f$ replaces the rough frequency cutoff $\chi_{B(0,R)}$ in the definition of $S_R^{\mathrm{ball}}f$ with the smoothed out version $m^{\alpha}_R(\xi) := \chi_{B(0,R)}(\xi)(1 - |\xi|^2/R^2)^{\alpha}$. The graph of this function is sketched in Figure~\ref{fig: Bochner Riesz}. For $\alpha > 0$ small, $m^{\alpha}_R$ is continuous, but fails to be differentiable at the boundary of $B(0,R)$. As $\alpha$ increases, it becomes more and more regular along the boundary of the ball.

The Bochner--Riesz means provide a natural mollification of the ball summation method. Given $1 \leq p \leq \infty$, this prompts the question of whether
\begin{equation}\label{eq: BR}
    \|B_R^{\alpha}f - f\|_{L^p(\R^n)} \to 0 \qquad \textrm{as $R \to \infty$} \qquad \textrm{holds for all $f \in L^p(\R^n)$.}
\end{equation}
As for the ball summation method, Plancherel's theorem shows \eqref{eq: BR} is trivially true for $p=2$. The \textit{Bochner--Riesz Conjecture} asserts that, as soon as $\alpha > 0$, there is a range\footnote{The precise statement of the conjecture posits that this should hold for $\frac{2n}{n+1} < p < \frac{2n}{n-1}$, but we shall not concern ourselves with the exact numerology here.} of exponents $p$ for which \eqref{eq: BR} holds, beyond the trivial $p = 2$ case.

%%%%%%%%%%%%%%%%%%%%%%%%%%%%%%%%%%%%%%%%%%%%%%%%%%%%%%%%%%%%%%%%%%%%%%%%%%%%%%%%%%%%%%%%%%%%%%%%

%    A counterexample to the Bochner--Riesz Conjecture

%%%%%%%%%%%%%%%%%%%%%%%%%%%%%%%%%%%%%%%%%%%%%%%%%%%%%%%%%%%%%%%%%%%%%%%%%%%%%%%%%%%%%%%%%%%%%%%%

\subsection{A counterexample to the Bochner--Riesz Conjecture?}  One lesson of Fefferman's work is that it is healthy to have a good dose of scepticism. Rather than try to prove a positive result on the Bochner--Riesz Conjecture, we could first try to construct a counterexample. 

As it stands, Fefferman's argument, based on the Perron tree construction of a Kakeya set of measure zero, is too weak to contradict the Bochner--Riesz Conjecture. Any tiny additional degree of smoothing (that is, any choice of $\alpha > 0$) is enough to dampen and overcome the pile up effect described in \S\ref{sec: Fefferman}. But what would happen if there was a way to construct even smaller Kakeya sets? 

Despite having zero measure, the planar Kakeya set arising from Perron tree construction is still relatively large: it has \textit{dimension} $2$. Here \textit{dimension} can either refer to the \textit{Minkowski} or the \textit{Hausdorff dimension}. We shall (briefly and informally) review these concepts in \S\ref{subsec: WZ proof} below. For our present purposes, however, it suffices to intuitively think of a $2$-dimensional subset of $\R^2$ as large and filling out space (whilst still possibly having zero measure). 

Suppose there exists a Kakeya set $K \subseteq \R^n$ which not only has measure zero, but is lower dimensional; say the Minkowski dimension of $K$ is less than $n$. This pathological object would pack all the line segments into a truly tiny space. We could then apply Fefferman's argument to \textit{this} Kakeya set. The pile up effect for our new $K$ would be monstrous: far worse than the Perron tree construction. In particular, it would overwhelm the effect of smoothing, at least for small values of $\alpha > 0$, and give a counterexample to the Bochner--Riesz Conjecture. 

%%%%%%%%%%%%%%%%%%%%%%%%%%%%%%%%%%%%%%%%%%%%%%%%%%%%%%%%%%%%%%%%%%%%%%%%%%%%%%%%%%%%%%%%%%%%%%%%

%    Dimensions of Kakeya sets

%%%%%%%%%%%%%%%%%%%%%%%%%%%%%%%%%%%%%%%%%%%%%%%%%%%%%%%%%%%%%%%%%%%%%%%%%%%%%%%%%%%%%%%%%%%%%%%%

\section{Dimensions of Kakeya sets}

%%%%%%%%%%%%%%%%%%%%%%%%%%%%%%%%%%%%%%%%%%%%%%%%%%%%%%%%%%%%%%%%%%%%%%%%%%%%%%%%%%%%%%%%%%%%%%%%

%    The Kakeya Conjecture

%%%%%%%%%%%%%%%%%%%%%%%%%%%%%%%%%%%%%%%%%%%%%%%%%%%%%%%%%%%%%%%%%%%%%%%%%%%%%%%%%%%%%%%%%%%%%%%%

\subsection{The Kakeya Conjecture} The existence of low dimensional Kakeya sets would have profound consequences for harmonic analysis, providing counterexamples to many important problems, including the Bochner--Riesz Conjecture discussed above. However, it seems very difficult to construct a Kakeya set which is significantly smaller than that coming from the Perron tree construction. This led to the widely held belief that low dimensional Kakeya sets do not exist.

\begin{conjecture}\label{conj: Kakeya 1} Any Kakeya set $K \subseteq \R^n$ has dimension $n$. 
\end{conjecture}

Here, as before, \textit{dimension} refers either to the \textit{Minkowski} or the \textit{Hausdorff} dimension. At its heart the meaning of the conjecture is simple: although Kakeya sets can have measure zero, they cannot be \textit{too} small. 

If the Kakeya Conjecture failed, then there would exist a Kakeya set in $\R^n$ of small dimension (less than $n$). Following the discussion of the previous section, this set could be used to produce a Fefferman-style counterexample to the Bochner--Riesz Conjecture. In this way, we see that
\begin{equation}\label{eq: BR implies Kak}
    \textrm{Bochner--Riesz Conjecture} \qquad \implies \qquad \textrm{Kakeya Conjecture.} 
\end{equation}
Many other central conjectures in Fourier analysis are known to imply the Kakeya Conjecture, following a similar scheme. But the significance of the Kakeya Conjecture stretches far beyond Fourier analysis! There are important problems in PDE, combinatorics and number theory which are also known to either imply the Kakeya Conjecture or some variant. 

%We present an array of such implications and connections (both formal and moral) in Figure~\ref{fig: connections}. 

%\begin{figure}
%    \centering
%    \includegraphics[width=\linewidth]{Images_PDF/connections2.png}
%    \caption{The web of problems related to the Kakeya Conjecture.}
 %   \label{fig: connections}
%\end{figure}

The Kakeya Conjecture in the plane is not too difficult. In the 1970s, Davies \cite{Davies1971} and later C\'ordoba \cite{Cordoba1977} gave two very different proofs showing that any planar Kakeya set has dimension $2$.\footnote{See, for instance, \cite{Wolff1999} and \cite{Mattila} for expositions of these results.} Correspondingly, the Bochner--Riesz Conjecture is also known to be true in two dimensions, and in general we have a relatively good understanding of Fourier analysis in the planar setting. By contrast, the Kakeya Conjecture in three (and higher) dimensional space is much, much harder. The additional directions in three dimensions complicate matters enormously. Correspondingly, our understanding of many central questions in Fourier analysis in $\R^3$ is severely lacking. The Kakeya Conjecture in higher dimensions has been a major bottleneck in harmonic analysis over the last 50 years. For this reason, it is considered one of the most important open problems in the area.

%The Bochner--Riesz Conjecture is known for $n = 2$, but remains a major open problem in dimensions $n \geq 3$. 

%%%%%%%%%%%%%%%%%%%%%%%%%%%%%%%%%%%%%%%%%%%%%%%%%%%%%%%%%%%%%%%%%%%%%%%%%%%%%%%%%%%%%%%%%%%%%%%%

%    Spectacular progress

%%%%%%%%%%%%%%%%%%%%%%%%%%%%%%%%%%%%%%%%%%%%%%%%%%%%%%%%%%%%%%%%%%%%%%%%%%%%%%%%%%%%%%%%%%%%%%%%

\subsection{Spectacular progress} Many mathematicians have studied the Kakeya Conjecture in three (and also higher) dimensions over the last 50 years. By their efforts, deep connections were discovered between the problem and questions in PDE, arithmetic combinatorics, number theory, and beyond \cite{KT2002, Wolff1999}. Partial results on the Kakeya Conjecture have led to dramatic breakthroughs on other problems.\footnote{To highlight one particularly striking instance, the \textit{multilinear} Kakeya theory of Bennett--Carbery--Tao~\cite{BCT2006} led to the resolution of the longstanding main conjecture in Vinogradov's mean value theorem~\cite{BDG2016} in analytic number theory.} However, despite remarkable successes, a full solution, even for $n=3$, always appeared far out of reach. The Kakeya Conjecture achieved a certain notoriety and became a byword for a difficult, near impossible challenge. This was the state of affairs until February 2025, when Hong Wang and Joshua Zahl \cite{WZ} announced the final piece of a spectacular full solution of the three dimensional problem! 

\begin{theorem}[Wang--Zahl \cite{WZsticky, WZ}]\label{thm: WZ} Any Kakeya set $K \subseteq \R^3$ has dimension $3$. 
\end{theorem}

The monumental proof of Theorem~\ref{thm: WZ} is hundreds of pages long and is spread over multiple papers. It combines a rich array of ideas, and relies on contributions from many different mathematicians over the last 30 years. However, their arguments also incorporate a wealth of new and deep insights and mark a major leap in our understanding. The work of Wang and Zahl is a great and historic achievement in modern mathematics. 

%Researchers are now working hard to fully digest Wang and Zahl's innovations. It is likely that their work will lead to breakthroughs on a wide range important problems, such as the Bochner--Riesz Conjecture, and the rich array of other questions appearing in Figure~\ref{fig: connections}. Their work also leaves open the Kakeya Conjecture in $\R^n$ for $n \geq 4$. Significant new features appear in the four (and higher) dimensional Kakeya problem, and the hunt is on to understand the extent to which the Wang--Zahl approach can be extended. The next chapter of this story will be an exciting one!

%%%%%%%%%%%%%%%%%%%%%%%%%%%%%%%%%%%%%%%%%%%%%%%%%%%%%%%%%%%%%%%%%%%%%%%%%%%%%%%%%%%%%%%%%%%%%%%%

%    Elements of the proof

%%%%%%%%%%%%%%%%%%%%%%%%%%%%%%%%%%%%%%%%%%%%%%%%%%%%%%%%%%%%%%%%%%%%%%%%%%%%%%%%%%%%%%%%%%%%%%%%

\subsection{Elements of the proof}\label{subsec: WZ proof} We conclude with a discussion of some of the features of the Wang--Zahl proof. Given the length and complexity of their arguments, we only provide the roughest of sketches. 

%%%%%%%%%%%%%%%%%%%%%%%%%%%%%%%%%%%%%%%%%%%%%%%%%%%%%%%%%%%%%%%%%%%%%%%%%%%%%%%%%%%%%%%%%%%%%%%%

%    Discretisation

%%%%%%%%%%%%%%%%%%%%%%%%%%%%%%%%%%%%%%%%%%%%%%%%%%%%%%%%%%%%%%%%%%%%%%%%%%%%%%%%%%%%%%%%%%%%%%%%

\medskip \noindent \textit{Discretisation.} The first step is a standard reformulation of the problem. Suppose $K \subseteq \R^3$ is a Kakeya set and, for $0 < \delta < 1$, let $N_{\delta} K$ denote the $\delta$-neighbourhood of $K$. We may assume $K$ has measure zero (since otherwise it automatically has dimension $3$). It then follows that
\begin{equation}\label{eq: discrete 1}
    |N_{\delta} K| \to 0 \qquad \textrm{as $\delta \to 0_+$,}
\end{equation}
where here the notation $| \,\cdot\,|$ is used to denote $3$-dimensional volume. 

The Minkowski dimension quantifies the rate of convergence in \eqref{eq: discrete 1}. In order to prove our Kakeya set $K$ has Minkowski dimension $3$, given $\varepsilon > 0$, we need to show that there exists some constant $c_{\varepsilon} > 0$ such that 
\begin{equation}\label{eq: discrete 2}
      |N_{\delta} K| \geq c_{\varepsilon} \delta^{\varepsilon} \qquad \textrm{for all $\delta > 0$.}
\end{equation}
In particular, although $|N_{\delta} K| \to 0$ as $\delta \to 0_+$, the rate of convergence is much slower than any power of $\delta$. In this sense, $K$ is `large'. We remark that Hausdorff dimension provides a more nuanced notion of size, which relies on analysing $K$ simultaneously at many scales, rather than just one fixed scale $\delta$. For the purposes of this sketch, however, we stick to the simpler notion of Minkowski dimension.

When we form the $\delta$-neighbourhood $N_{\delta} K$ of $K$, each line segment in $K$ fattens up into a solid shape which is roughly a cylinder of length $1$ and radius $\delta$. We call these cylinders \textit{$\delta$-tubes}. We expect many of the $\delta$-tubes formed from lines in $K$ to essentially coincide. For this reason, we sample from $K$ a finite family of line segments $\cL$ with directions forming a maximal, $\delta$-separated set. The corresponding $\delta$-tubes are then \textit{essentially distinct} (they differ significantly from one another). To succinctly describe this setup, we make the following definition.

\begin{definition} For $0 < \delta < 1$, we say that a set $\T$ of $\delta$-tubes is \textit{direction-separated} if the directions of the core lines of the tubes in $\T$ form a $\delta$-separated set.
\end{definition} 
In light of the above, the proof of \eqref{eq: discrete 2} reduces to showing the following. 

\begin{theorem}[Kakeya Conjecture, discretised version \cite{WZsticky, WZ}]\label{thm: discrete Kak} Given any $\varepsilon > 0$, there exists some $c_{\varepsilon} > 0$ such that for all $0 < \delta < 1$ the inequality
\begin{equation}\label{eq: discrete 3}
    \Big|\bigcup_{T \in \T} T\Big| \geq c_{\varepsilon} \delta^{\varepsilon} \sum_{T \in \T} |T|
\end{equation}
holds whenever $\T$ is a direction-separated set of $\delta$-tubes in $\R^3$.
\end{theorem} 

If the tubes in $\T$ were disjoint, then we would have $\big|\bigcup_{T \in \T} T\big| =\sum_{T \in \T} |T|$. Thus, \eqref{eq: discrete 3} can be interpreted as saying the tubes in $\T$ are `almost' disjoint. 

To pass from \eqref{eq: discrete 3} to \eqref{eq: discrete 2}, we simply take $\T$ to be the set of tubes formed from the line segments belonging to $\cL$ as defined earlier. On the one hand, $N_{\delta} K \supseteq \bigcup_{T \in \T} T$. On the other hand, maximality ensures that $\sum_{T \in \T} |T| \sim 1$. These observations combine with \eqref{eq: discrete 3} to give \eqref{eq: discrete 2}. The advantage of the reformulation in Theorem~\ref{thm: discrete Kak} is that, rather than studying a \textit{continuum} of line segments in $K$, here we study a \textit{finite} set of $\delta$-tubes. This is a form of \textit{discretisation}, a common tool in Fourier analysis and geometric measure theory.\footnote{We also saw an instance of discretisation in the construction of the Fefferman counterexample. }
%%%%%%%%%%%%%%%%%%%%%%%%%%%%%%%%%%%%%%%%%%%%%%%%%%%%%%%%%%%%%%%%%%%%%%%%%%%%%%%%%%%%%%%%%%%%%%%%

%    Exploiting distinct directions

%%%%%%%%%%%%%%%%%%%%%%%%%%%%%%%%%%%%%%%%%%%%%%%%%%%%%%%%%%%%%%%%%%%%%%%%%%%%%%%%%%%%%%%%%%%%%%%%

\medskip \noindent \textit{Exploiting distinct directions.} It is surprisingly difficult to directly leverage the \textit{direction-separated} hypothesis from the Kakeya set definition. In practice, one often works with weaker, and more malleable, conditions. We noted earlier that direction-separated tubes are \textit{essentially distinct}. This is formally defined as follows. 

\begin{definition} Let $0 < \delta < 1$ and $\T$ be a set of $\delta$-tubes. We say the tubes in $\T$ are \textit{essentially distinct} if $|T_1 \cap T_2| \leq |T_1|/2$ for all $T_1$, $T_2 \in \T$ with $T_1 \neq T_2$.     
\end{definition}

It is natural to ask whether the conclusion of Theorem~\ref{thm: discrete Kak} continues to hold under the weaker hypothesis that the tubes in $\T$ are essentially distinct. However, there is an easy counterexample. Assume $0< \delta < 1$ satisfies $\delta^{-1} \in \N$ and let 
\begin{equation*}
A := \{(\delta j, 0, 0) : 1 \leq j \leq \delta^{-1} \} \quad \textrm{and} \quad B := \{(\delta j, 1, 0) : 1 \leq j \leq \delta^{-1} \}
\end{equation*}
 be sets of equally spaced points along a pair of parallel lines in $S := [0,1]^2 \times \{0\}$. Form a family $\cL$ of $\delta^{-2}$ line segments by joining every point in $A$ to a point in $B$. If $\T$ denotes the corresponding set of $\delta$-tubes, then $\#\T = \delta^{-2}$ and the tubes in $\T$ are essentially distinct. However, $T \subseteq N_{\delta} S$ for all $T \in \T$ and, moreover, $|\bigcup_{T \in \T} T| \sim \delta$. Thus, \eqref{eq: discrete 3} fails in this case. 

In the above example, the $\delta$-tubes concentrate in $N_{\delta} S$, which is essentially a $1 \times 1 \times \delta$ rectangular prism. What happens if we rule out this behaviour?

\begin{definition}[\cite{Wolff1995}] A set $\T$ of $\delta$-tubes in $\R^3$ satisfies the \textit{Wolff axiom} if 
\begin{equation}\label{eq: Wolff}
    \#\{T \in \T : T \subseteq R \} \leq \delta^{-2}|R|
\end{equation}
holds whenever $R \subseteq \R^3$ is a rectangular prism. 
\end{definition}

If $\T$ is direction-separated, then it is not difficult to see that $\T$ automatically satisfies the Wolff axiom. It is natural to ask whether the Wolff axiom is sufficient to ensure the conclusion of Theorem~\ref{thm: discrete Kak}. It turns out that this is indeed the case. 

\begin{theorem}[Strong Kakeya Conjecture \cite{WZsticky, WZ2025, WZ}]\label{thm: strong Kak} Given any $\varepsilon > 0$, there exists some $c_{\varepsilon} > 0$ such that for all $0 < \delta < 1$ the inequality
\begin{equation}\label{eq: strong Kak}
    \Big|\bigcup_{T \in \T} T\Big| \geq c_{\varepsilon} \delta^{\varepsilon}\sum_{T \in \T} |T|
\end{equation}
holds whenever $\T$ is a collection of $\delta$-tubes in $\R^3$ which satisfy the Wolff axiom.
\end{theorem} 

At this point we encounter the first major indication of the subtlety of the Kakeya problem: Theorem~\ref{thm: strong Kak} \textit{fails} over the complex field! In \cite[\S13]{KLT2000}, authors considered the compact piece of the \textit{Heisenberg group} in $\C^3$, defined by
\begin{equation*}
    H := \big\{(z_1, z_2, z_3) \in \C^3 : \fkI(z_1) = \fkI(z_2\overline{z}_3),\, |z_1|, |z_2|, |z_3| \leq 2 \big\}.
\end{equation*}
It is easy to see that $H$ contains a $4$ (real) parameter family of complex line segments
\begin{equation*}
    \ell_{a, b, w} := \{(z, w + az, z\overline{w} + b) : z \in \C,\, |z| \leq 1/2\},
\end{equation*} 
where $a$, $b \in \R$ and $w \in \C$ with $|a|$, $|b|$, $|w| \leq 1$. From this, one may construct a set $\T$ of complex $\delta$-tubes with $\#\T \sim \delta^{-4}$ which satisfies the natural complex analogue of the Wolff axiom. Since each complex $\delta$-tube $T$ has measure $|T| \sim \delta^4$, it follows that $\sum_{T \in \T} |T| \sim 1$. However, $T \subseteq N_{\delta} H$ for all $T \in \T$ and, morover, $\big|\bigcup_{T \in \T} T\big| \sim \delta$. Thus, $\T$ fails \eqref{eq: strong Kak}.

In light of the above example, any proof of Theorem~\ref{thm: strong Kak} must somehow distinguish between the real and complex field. The definition of the Heisenberg group relies on the notion of the real and imaginary parts of a complex number, and the accompanying conjugation operator $z \mapsto \bar{z}$. These definitions in turn rest on the existence of the half-dimensional subfield $\R \subseteq \C$. There is no apparent analogue of this structure in the real line. Indeed, a famous result of Edgar--Miller \cite{EM2003} shows that any Borel subring of $\R$ either has Hausdorff dimension $0$ or is the whole of $\R$. This fact will play an important role later in the proof.

%%%%%%%%%%%%%%%%%%%%%%%%%%%%%%%%%%%%%%%%%%%%%%%%%%%%%%%%%%%%%%%%%%%%%%%%%%%%%%%%%%%%%%%%%%%%%%%%

%    Reducing to the highly symmetric case

%%%%%%%%%%%%%%%%%%%%%%%%%%%%%%%%%%%%%%%%%%%%%%%%%%%%%%%%%%%%%%%%%%%%%%%%%%%%%%%%%%%%%%%%%%%%%%%%

\medskip \noindent \textit{Reducing to a symmetric case.} To prove Theorem~\ref{thm: strong Kak}, we can attempt to isolate configurations of tubes which are extremal, or optimal, in the sense that $\big|\bigcup_{T \in \T} T\big|$ is minimal over all sets of tubes $\T$ of a fixed cardinality satisfying the Wolff axiom. If we can prove that the Strong Kakeya Conjecture holds for such extremal configurations, it will automatically hold in all cases. 

It is natural to expect an optimal configuration of tubes should saturate the Wolff axiom hypothesis of Theorem~\ref{thm: strong Kak}. Indeed, if there was any leeway in our hypotheses, we could hope to exploit it to construct a tighter packing of tubes. This leads to the definition of a \textit{sticky} set of tubes.  

\begin{definition}[Sticky, informal] For each intermediate scale $\delta \leq \rho \leq 1$, we consider fattening the $\delta$-tubes $T$ into $\rho$-tubes $T_{\rho}$. We expect many of the $\rho$-tubes formed in this way to essentially coincide. We therefore extract a maximal family of essentially distinct $\rho$-tubes $\T_{\rho}$. We say the collection of $\delta$-tubes is \textit{sticky} if 
\begin{equation}\label{eq: sticky}
    \#\{T \in \T : T \subseteq T_{\rho}\} \sim (\rho/\delta)^2 \qquad \textrm{for all $T_{\rho} \in \T_{\rho}$ and all $\delta \leq \rho \leq 1$.}
\end{equation}
In other words, each $\rho$-tube in $\T_{\rho}$ contains the maximal number of $\delta$-tubes permitted by the Wolff axiom condition \eqref{eq: Wolff}. 
\end{definition}

Sticky sets of tubes are highly symmetric from the perspective of metric dimension. The condition \eqref{eq: sticky} relates the distribution of the tubes $\T$ at different scales $\delta \leq \rho \leq 1$ and is an example of \textit{statistical self-similarity}: while the set is not literally self-similar (in the sense of the middle third Cantor set, for example), certain quantitative properties of the set are identical at many different scales. Our experience and intuition suggest optimal configurations should be highly symmetric. Moreover, the Kakeya Conjecture is both a quantitative and a multi-scale problem: we are interested in the size, and not the exact shape, of Kakeya sets and we need to prove results \textit{for all} scales $\delta$. The relevant notions of symmetry should reflect these fundamental features. This inexorably leads to statistical self-similarity. Thus, sticky sets of tubes arise as a natural and important special case in the study of the Kakeya Conjecture.

\begin{theorem}[Sticky reduction \cite{WZ}]\label{thm: stricky red} In order to prove Theorem~\ref{thm: strong Kak}, it suffices to assume, in addition to the existing hypotheses, that $\T$ is sticky.
\end{theorem}

This is the most significant innovation in the work of Wang--Zahl. Loosely speaking, Theorem~\ref{thm: stricky red} formalises our intuition that optimal configurations of tubes for Theorem~\ref{thm: strong Kak} correspond to the highly symmetric arrangements for which the Wolff axiom is saturated. However, making this intuition precise is monumentally challenging. Reduction to the sticky case was previously explored by Wolff and Katz--\L aba--Tao \cite{KLT2000} in the early 2000s, but with limited success. Since then very little progress had been made on this problem. Indeed, prior to Wang--Zahl~\cite{WZ}, there was no hint of a viable mechanism for reducing to the sticky case; to the contrary, evidence was mounting \cite{KZ2019} that such a reduction was infeasible. The proof of Theorem~\ref{thm: stricky red} forms the entirety of the 127 page preprint \cite{WZ} and introduces very important new ideas. 

It is unfortunately outside the modest scope of this article to delve into the proof of Theorem~\ref{thm: stricky red}. However, we mention that it relies on an intricate induction or self-improving mechanism. For this reason, it is essential to work with the Wolff axiom formulation of the problem as in the Strong Kakeya Conjecture. Indeed, the hypotheses must be sufficiently general in order to close the inductive arguments.

%%%%%%%%%%%%%%%%%%%%%%%%%%%%%%%%%%%%%%%%%%%%%%%%%%%%%%%%%%%%%%%%%%%%%%%%%%%%%%%%%%%%%%%%%%%%%%%%

%    A complex variant

%%%%%%%%%%%%%%%%%%%%%%%%%%%%%%%%%%%%%%%%%%%%%%%%%%%%%%%%%%%%%%%%%%%%%%%%%%%%%%%%%%%%%%%%%%%%%%%%

%\medskip \noindent \textit{A complex variant} The next key observation is that the Strong Kakeya Conjecture from Theorem~\ref{thm: strong Kak} \textbf{fails} in the setting of the complex field $\C^3$.  \jh{go through this!!} 

%The existence of the Heisenberg group example implies that any proof of Theorem~\ref{thm: strong Kak} must somehow distinguish between the real field $\R$ and the complex field $\C$.

%%%%%%%%%%%%%%%%%%%%%%%%%%%%%%%%%%%%%%%%%%%%%%%%%%%%%%%%%%%%%%%%%%%%%%%%%%%%%%%%%%%%%%%%%%%%%%%%

%    Forbidden symmetries

%%%%%%%%%%%%%%%%%%%%%%%%%%%%%%%%%%%%%%%%%%%%%%%%%%%%%%%%%%%%%%%%%%%%%%%%%%%%%%%%%%%%%%%%%%%%%%%%

\medskip \noindent \textit{Forbidden symmetries.} It remains to prove the Strong Kakeya Conjecture under the additional hypothesis that the set of tubes $\T$ is sticky. To treat this case, Wang--Zahl \cite{WZsticky, WZ2025} were inspired by an approach originally proposed by Katz and Tao, and shared with the public in 2014 in a talk and blogpost \cite{Tao_blog} by Tao. At that time, Katz--Tao were unable to rigorously implement their proposal. Indeed, the implementation in Wang--Zahl \cite{WZsticky} involves substantial new ideas and the use of mathematical tools and techniques which were unavailable back in 2014. 

Katz--Tao proposed a proof by contradiction. We assume there exists a sticky set of tubes $\T$ which satisfies the hypotheses of the Strong Kakeya Conjecture, but for which $K := \bigcup_{T \in \T} T$ has very small volume. Moreover, we may assume the volume of $K$ is minimal over all such sets of tubes $\T$ of a fixed cardinality. 

The stickiness hypothesis means that $K$ is highly symmetric. Under the minimal volume condition, the symmetries are further amplified. In particular, following landmark ideas introduced by Katz--\L aba--Tao~\cite{KLT2000}, one can show the tubes in $\T$ must have extremely rich structural properties. Katz--Tao \cite{Tao_blog} proposed using these structural properties to encode highly symmetric \textit{algebraic} structures inside the set $K$. This is in part motivated by the earlier observation that any proof of Theorem~\ref{thm: strong Kak} must distinguish between the real and complex field (and, consequently, it is natural to expect algebra to play a role at some level).

One of the most symmetric objects in mathematics is a commutative ring: it simultaneously admits \textit{two} symmetry groups! By analysing the interplay between the structure of the Kakeya set at two different scales, Katz--Tao \cite{Tao_blog} sketched a method for encoding a ring-like object $R$ into $K$. More precisely, they proposed encoding an \textit{approximate subring} $R$ of $\R$ into $K$. This is a set which is approximately closed under the addition and multiplication operations, in the sense that sum set $R + R := \{x + y : x, y \in R\}$ and product set $R \cdot R := \{x \cdot y : x, y \in R\}$ have size comparable to that of $R$. Furthermore, this approximate subring would have Hausdorff dimension $0 < \dim_{\mathrm{H}} R < 1$. 

At this point, the proof sketch is nearing its endgame. We have already mentioned that the Edgar--Miller theorem \cite{EM2003} rules out the existence of subrings of intermediate Hausdorff dimension. This would give the desired contradiction, except we are dealing with \textit{approximate} subrings rather than \textit{bona fide} subrings. However, there is a more quantitative version of the Edgar--Miller theorem called \textit{Bourgain's sum-product theorem} \cite{Bourgain2003, Bourgain2010}, which does apply in the approximate setting. Using this, Katz--Tao hoped to obtain a contradiction and thereby conclude the proof of the Strong Kakeya Conjecture in the sticky case.

As mentioned earlier, the rigorous proof in~\cite{WZsticky} differs significantly from that originally envisioned by Katz--Tao \cite{Tao_blog}. Most notably, the connection with approximate subrings is far less direct. Rather than appeal to Bourgain's sum-product theorem \cite{Bourgain2010}, Wang--Zahl~\cite{WZsticky} invoke various recent results in \textit{projection theory}, a subfield of geometric measure theory, to arrive at a contradiction. These tools ultimately do rely on Bourgain's sum-product theorem \cite{Bourgain2010} but, taken as a whole, the scheme of the proof is far more complex than that sketched in \cite{Tao_blog}.

%%%%%%%%%%%%%%%% BIBLIOGRAPHY %%%%%%%%%%%%%%%%%%%%%%%%%%%%%%%%%%%%%%%%%%%%%%%%%%%%%%%%%%%%%%%%%%

\bibliography{Reference}
%%%
\end{document}